\documentclass[11pt]{article}
\usepackage{amssymb}
\usepackage{latexsym}
\input amssym.def
\input amssym
\newtheorem{theorem}{Theorem}[section]
\newtheorem{lemma}[theorem]{Lemma}

\newtheorem{example}[theorem]{Example}
\newtheorem{proposition}[theorem]{Proposition}
\newtheorem{corollary}[theorem]{Corollary}

\newtheorem{remark}[theorem]{Remark}

\newcommand{\halmos}{\rule{1ex}{1.4ex}}

\newcommand{\bea}{\begin{eqnarray}}
\newcommand{\epf}{\hspace*{\fill}\mbox{$\halmos$}}
\newcommand{\eea}{\end{eqnarray}}
\newcommand{\nn}{\nonumber \\}
\newcommand{\be}{\begin{equation}}

\newcommand{\ee}{\end{equation}}

\usepackage{amssymb}
\title{{The Rogers--Selberg recursions, the Gordon--Andrews identities
and intertwining operators} }
\author{S. Capparelli \footnote{S.C. gratefully acknowledges
partial support {}from MIUR (Ministero dell'Istruzione, dell'Universit\`a
e della Ricerca).}, J. Lepowsky\footnote{J.L. and A.M.  gratefully
acknowledge partial support {}from NSF grant DMS--0070800.} and A. Milas
\footnotemark[\value{footnote}] }
\begin{document}
\date{}
\bibliographystyle{alpha}
\maketitle
\begin{abstract}
Using the theory of intertwining operators for vertex operator
algebras we show that the graded dimensions of the principal subspaces
associated to the standard modules for $\widehat{\goth{sl}(2)}$
satisfy certain classical recursion formulas of Rogers and Selberg.
These recursions were exploited by Andrews in connection with Gordon's
generalization of the Rogers--Ramanujan identities and with Andrews'
related identities.  The present work generalizes the authors'
previous work on intertwining operators and the Rogers--Ramanujan
recursion.
\end{abstract}
\maketitle

\renewcommand{\theequation}{\thesection.\arabic{equation}}
\renewcommand{\thetheorem}{\thesection.\arabic{theorem}}
\setcounter{equation}{0}
\setcounter{theorem}{0}
\setcounter{section}{0}

\section{Introduction}
This paper is a continuation of \cite{CLM}, to which we refer the
reader for background and for our motivation.

We continue our study of a relationship between intertwining
operators, in the sense of vertex operator algebra theory, associated
to standard (integrable highest weight) modules for
$\widehat{\goth{sl}(2)}$, and the corresponding principal subspaces
(\cite{FS1}--\cite{FS2}).  Let us recall our main result in
\cite{CLM}.

Consider the principal subspaces $W(\Lambda_0)$ and $W(\Lambda_1)$
associated to the level $1$ standard $\widehat{\goth{sl}(2)}$--modules
$L(\Lambda_0)$ and $L(\Lambda_1)$, respectively (see Section 2 below).
All these spaces are (doubly) graded by the eigenvalues of the
operator $\alpha/2$, $\alpha$ being the positive root of
$\goth{sl}(2)$, and the operator $L(0)$ coming {}from the Virasoro
algebra.  The corresponding graded dimensions (the
``characters''---the generating functions of the dimensions of the
homogeneous subspaces) of $W(\Lambda_0)$ and $W(\Lambda_1)$ are then
given by:
\bea {\chi}_{W(\Lambda_0)}(x,q) &=&{\rm tr}_{W(\Lambda_0)}
x^{\alpha/2} q^{L(0)},\nn
{\chi}_{W(\Lambda_1)}(x,q)&=& {\rm
tr}_{W(\Lambda_1)} x^{\alpha/2} q^{L(0)}; \nonumber
\eea
the variables
$x$ and $q$ are formal.  It is convenient to use the slightly modified
graded dimension
$${\chi}'_{W(\Lambda_1)}(x,q)=
x^{-1/2}q^{-1/4}{\rm tr}_{W(\Lambda_1)} x^{\alpha/2} q^{L(0)},$$
and for cosmetic reasons we correspondingly let
$${\chi}'_{W(\Lambda_0)}(x,q)={\chi}_{W(\Lambda_0)}(x,q).$$
We proved \cite{CLM} that
these graded dimensions satisfy the linear $q$--difference equation
\be
\label{rrrecursion}
{\chi}'_{W(\Lambda_0)}(x,q)-{\chi}'_{W(\Lambda_1)}(x,q)=
xq{\chi}'_{W(\Lambda_1)}(xq,q).
\ee
We did this by constructing a certain exact sequence by means of
intertwining operators associated to irreducible modules for
$L(\Lambda_0)$ viewed as a vertex operator algebra.  Formula
(\ref{rrrecursion}) together with the formula
\be
{\chi}'_{W(\Lambda_1)}(x,q)={\chi}'_{W(\Lambda_0)}(xq,q)
\ee
yields
\be \label{rrr}
{\chi}'_{W(\Lambda_0)}(x,q)-{\chi}'_{W(\Lambda_0)}(xq,q)=
xq{\chi}'_{W(\Lambda_0)}(xq^2,q), \ee
which is the classical Rogers--Ramanujan recursion (cf. \cite{A}),
asserted for the graded dimension of $W(\Lambda_0)$.  In particular,
this easily implies (cf. \cite{A}) that
\be \label{solnofrrr}
{\chi}'_{W(\Lambda_0)}(x,q)=\sum_{n \geq 0} \frac{x^n q^{n^2}}{(q)_n},
\ee
where $(q)_n=(1-q) \cdots (1-q^n)$.  The specializations $x=1$ and
$x=q$ of (\ref{solnofrrr}) give the sum sides of the two
Rogers--Ramunujan identities (cf. \cite{A}).

Gordon's identities \cite{G} generalize the Rogers--Ramanujan
identities.  They can be stated in the following form:
\begin{theorem} \label{gi}
Let $l \geq 2$ and let $1 \leq t \leq l$.  The number of partitions of
a nonnegative integer $n$ into parts not congruent to $0$, $\pm t$
${\rm mod} \ (2l+1)$ is equal to the number of partitions $n=b_1+
\cdots +b_s$, $b_1 \geq b_2 \geq \cdots \geq b_s > 0$, such that
$b_i-b_{i+l-1} \geq 2$ (the ``difference two at distance $l-1$''
condition) and at most $t-1$ of the $b_i$ equal $1$.
\end{theorem}

In \cite{A1} (cf. \cite{A}) Andrews found an elegant proof of Gordon's
identities by using certain $q$--hypergeometric--type series and a
family of recursions (\cite{RR}, \cite{S}) expressed as a system of
$q$--difference equations (\cite{S}, \cite{A1}).  We will call these
recursions, or equivalently, $q$--difference equations, the {\em
Rogers--Selberg recursions}.

It is easy to see that the generating function of the number of
partitions of a nonnegative integer into parts not congruent to $0$,
$\pm t$ ${\rm mod} \ (2l+1)$ can be expressed as
$$\prod_{i>0, \;\; i \neq 0,\pm t \ {\rm mod} \ (2l+1)}
\frac{1}{1-q^i},$$
the {\em product side} of Gordon's identities.  Andrews \cite{A4}
discovered an ``analytic form'' of (the sum sides of) these
identities, and he derived it {}from the same difference equations
mentioned above.  Actually, this ``analytic form'' is a ``multisum
form,'' and Andrews exploited a refined version of the $q$--generating
function, involving the variable $x$, giving a generating function of
a refined version of the partitions counted in Gordon's form of the
identities.  We will continue to take the variables in these
``analytic'' expressions to be formal rather than complex.  With the
variable $x$ suitably specialized, these Gordon--Andrews identities
state (\cite{A4}; cf. \cite{A}):
\begin{theorem} \label{agi}
With $l$ and $t$ as in Theorem \ref{gi},
\bea \label{angordon}
&& \prod_{i > 0, \;\; i \neq 0, \pm t \ {\rm mod} \ (2l+1)}
\frac{1}{1-q^i} = \nn
&& = \sum_{N_1 \geq \cdots \geq N_{l-1} \geq 0}
\frac{q^{N_1^2 +
\cdots + N_{l-1}^2+N_t+N_{t+1}+\cdots + N_{l-1}}}{(q)_{N_1-N_2} \cdots
(q)_{N_{l-2}-N_{l-1}} (q)_{N_{l-1}}}.
\eea
\end{theorem}

We emphasize that for us, $q$ is a formal variable; in the literature,
one often takes $q$ to be a complex variable such that $|q|<1$.

A nonclassical approach to Gordon's identities was initiated and
developed by Lepowsky and Wilson in \cite{LW2}--\cite{LW4} by means of
``twisted vertex operators'' and ``$Z$--algebras''; they used this to
give a construction of bases for all the standard
$\widehat{\goth{sl}(2)}$--modules (or more particularly, bases for the
``vacuum spaces'' for the ``principal Heisenberg subalgebra''
\cite{LW1} of $\widehat{\goth{sl}(2)}$) compatible with Gordon's
identities in the case of the odd levels, and with certain
generalizations of Gordon's identities due to Andrews (\cite{A2},
\cite{A3}) and Bressoud (\cite{B1}, \cite{B2}) in the case of the even
levels.  Analogously, but with a somewhat different flavor, in
\cite{LP1} and \cite{LP2}, Lepowsky and Primc used untwisted vertex
operators and $Z$--algebras to construct bases for all the standard
$\widehat{\goth{sl}(2)}$--modules.  These twisted and untwisted
$Z$--algebra constructions (\cite{LW2}--\cite{LW4},
\cite{LP1}--\cite{LP2}) both led naturally to the
``difference--two--at--a--distance'' condition that had appeared in
Gordon's identities; in the untwisted case \cite{LP2} it led to a new
version of this condition.

Then Meurman and Primc \cite{MP0} used $Z$--algebras to complete a new
{\it proof} of Gordon's identities in the setting of
\cite{LW2}--\cite{LW4} (Lepowsky and Wilson had obtained a
$Z$--algebra proof of the identities in special cases, including the
Rogers--Ramunujan identities, in \cite{LW2}--\cite{LW3}).  On the
other hand, Feigin and Stoyanovsky (\cite{FS1}--\cite{FS2}) obtained
the sum sides of Gordon's identities as the graded dimensions of what
they called the ``principal subspaces'' associated to the standard
$\widehat{\goth{sl}(2)}$--modules of level $l-1$ (which we shall also
write as $k$, below).  In order to compute these ``characters,''
Feigin and Stoyanovsky considered ``quasi--particles,'' which are the
expansion coefficients of certain generalized ``fields'' (cf. Chapter
13 in \cite{DL}).  This approach was further developed by Georgiev in
\cite{Ge}. In addition, Feigin and Stoyanovsky interpreted the product
side of (\ref{angordon}) by using a geometric approach, via
infinite--dimensional analysis on a flag manifold. This approach was
further extended in \cite{FL} (see also \cite{FKLMM1}--\cite{FKLMM3}).

All these developments have concerned
$\widehat{\goth{sl}(2)}$--modules; there have been many other
developments in this spirit (see \cite{FF}, etc.).

In this paper we analyze the structure of the principal subspaces
associated to the standard modules for $\widehat{\goth{sl}(2)}$ by
using the theory of vertex operator algebras and intertwining
operators (\cite{FLM}, \cite{FHL}). Our main result is a construction
of $l$ exact sequences that yield a linear system of $q$--difference
equations equivalent to the Rogers--Selberg recursions (\cite{RR},
\cite{S}) mentioned above (see also \cite{A1}).  It is important to
say that, in order to obtain these exact sequences, we use certain
aspects of the theory of intertwining operators for vertex operator
algebras.  As in \cite{A1} and \cite{A4} (cf. \cite{A}), these
recursions in turn yield the sum sides of Gordon's identities (see
Theorem \ref{gi}) and of Andrews' identities (see Theorem \ref{agi}).

\renewcommand{\theequation}{\thesection.\arabic{equation}}
\renewcommand{\thetheorem}{\thesection.\arabic{theorem}}
\setcounter{equation}{0}
\setcounter{theorem}{0}

\section{Principal subspaces}

The setting and notation are as in \cite{CLM}.  Let
$\goth{g}=\goth{sl}(2, \mathbb{C})$ be the $3$--dimensional complex
simple Lie algebra with a standard basis $\{h, x_\alpha, x_{-\alpha}
\}$ and bracket relations
$$[h,x_{\alpha}]=2 x_{\alpha}, \ [h,x_{-\alpha}]=-2 x_{-\alpha}, \
[x_\alpha,x_{-\alpha}]=h.$$
We fix the Cartan subalgebra $\goth{h}=\mathbb{C}h$, which we identify
with its dual by means of the form $\langle x,y \rangle={\rm tr}(xy)$
for $x,y \in {{\goth g}}$. Take $\alpha \in {\goth{h}}$ to be the root
corresponding to the root vector $x_{\alpha}$, and take this root to
be positive; then $ \langle \alpha,\alpha \rangle=2$ and we have the
root space decomposition
$${\goth{g}}={\goth{n}}_- \oplus {\goth{h}} \oplus {\goth{n}}_+ ,$$
where ${\goth{n}}_\pm=\mathbb{C}x_{\pm \alpha}$.  Note that under our
identifications,
$$h = \alpha.$$
We shall use the affine Lie algebra
$$\widehat{\goth{sl}(2)}= \goth{sl}(2,\mathbb{C}) \otimes
\mathbb{C}[t,t^{-1}] \oplus \mathbb{C}c,$$
where $c$ is a nonzero central element and
$$[x \otimes t^m,y \otimes t^n]=[x,y] \otimes t^{m+n}+\langle x,y
\rangle m \delta_{m+n,0}c$$
for $x,y \in {\goth{g}}$ and $m,n \in {\mathbb Z}$.

For a dominant weight $\Lambda \in (\goth{h} \oplus \mathbb{C}c)^*$,
let $L(\Lambda)$ be the corresponding standard ($=$ integrable highest
weight) $\widehat{\goth{sl}(2)}$--module with highest weight
$\Lambda$ (cf. \cite{K}).  We say that an
$\widehat{\goth{sl}(2)}$--module has {\it level} $k \in \mathbb{C}$ if
$c$ acts as multiplication by $k$.  The standard module $L(\Lambda)$
(which is irreducible) has nonnegative integral level, given by
$\langle \Lambda,c \rangle$.
(The evaluation map $\langle \cdot,\cdot \rangle$ of $(\goth{h} \oplus
\mathbb{C}c)^*$ on $\goth{h} \oplus \mathbb{C}c$ extends the form
$\langle \cdot,\cdot \rangle$ on ${\goth{h}}={\goth{h}}^*$; also, for
$\mu, \nu \in (\goth{h} \oplus \mathbb{C}c)^*$, $\langle \mu,\nu
\rangle$ is defined to be $\langle{\mu}|_{\goth{h}},
\nu|_{\goth{h}}\rangle$.)
The highest weights of the level $1$
standard modules are the fundamental weights $\Lambda_0$ and
$\Lambda_1$, defined by: $\langle \Lambda_{i},c \rangle = 1$, $\langle
\Lambda_{i},h \rangle = \delta_{i,1}$ for $i=0,1$.  The highest
weights of the level $k$ standard modules are given by
$$k_0 \Lambda_0+k_1 \Lambda_1, \ \ \ k_0+k_1=k, $$
where $k_0$ and $k_1$ are nonnegative integers.

Throughout this paper we will write $x(n)$ for the action of $x
\otimes t^n \in \widehat{\goth{sl}(2)}$ on an
$\widehat{\goth{sl}(2)}$--module, for $x \in \goth{g}$ and $n \in
\mathbb{Z}$.

Let $P=\frac{1}{2}\mathbb{Z}\alpha$ be the weight lattice and
$Q=\mathbb{Z} \alpha$ the root lattice in $\goth{h}$.  Let
$\mathbb{C}[P]$ and $\mathbb{C}[Q]$ be the corresponding group
algebras, with bases $\{e^{\mu} \ | \ \mu \in P\}$ and $\{e^{\mu} \ |
\ \mu \in Q\}$.  Consider the subalgebra
$$\hat{\goth{h}}_{\mathbb Z}=\coprod_{m \in \mathbb{Z} \setminus \{0\}}
\goth{h} \otimes
t^m \oplus \mathbb{C}c$$
of $\widehat{\goth{sl}(2)}$, a Heisenberg algebra, meaning that its
commutator subalgebra is equal to its center, which is one--dimensional
(namely, $\mathbb{C}c$).  We will also need the subalgebra
\be
\label{hhat} \hat{\goth{h}}=\goth{h} \otimes \mathbb{C}[t,t^{-1}]
\oplus \mathbb{C}c.  \ee
We shall consider the $\hat{\goth{h}}$--module
$$M(1)=U( \hat{\goth{h}} ) \otimes_{U(\goth{h} \otimes
\mathbb{C}[t] \oplus \mathbb{C}c)} \mathbb{C},$$
where $\goth{h} \otimes \mathbb{C}[t]$ acts trivially on the
one--dimensional module $\mathbb{C}$ and $c$ acts as $1$.  It is well
known (\cite{FK}, \cite{Seg}; cf. \cite{FLM}) that
$$V_P=M(1) \otimes \mathbb{C}[P]$$
and its subspaces
$$V_Q=M(1) \otimes \mathbb{C}[Q]$$
and
$$V_{Q+\alpha/2}=M(1) \otimes e^{\alpha/2} \mathbb{C}[Q]$$
admit natural $\widehat{\goth{sl}(2)}$--module structure, via certain
vertex operators (recalled in \cite{CLM}), and that as
$\widehat{\goth{sl}(2)}$--modules,
$$V_{Q} \cong L(\Lambda_0) \ \ \mbox{and} \ \ V_{Q+\alpha/2} \cong
L(\Lambda_1).$$

It is much harder to obtain a similar construction for the higher
level standard $\widehat{\goth{sl}(2)}$--modules; this was done in
\cite{LP2}.  By tensoring $k$ level $1$ standard modules one obtains a
level $k$ module which is completely reducible and whose irreducible
components are level $k$ standard modules (cf. \cite{K}); what was
hard was to construct bases of such irreducible modules and to
determine their precise structure.  The action of $g \in
\widehat{\goth{sl}(2)}$ on a vector
$$v \in L(\Lambda_{i_1}) \otimes \cdots
\otimes L(\Lambda_{i_k}),$$
where each $i_r$ is either $0$ or $1$, is given by the usual comultiplication
\be \label{comult}
g \cdot v=\Delta(g)v=(g \otimes 1 \otimes \cdots \otimes 1+ \cdots
+ 1 \otimes \cdots  \otimes 1 \otimes g ) v,
\ee
and this action of course extends to $U(\widehat{\goth{sl}(2)})$.

Consider the subalgebra
\be
\label{nhat}
\goth{n}\widehat{}_+ =\goth{n}_+ \otimes \mathbb{C}[t,t^{-1}]
\ee
of $\widehat{\goth{sl}(2)}$.  The {\em principal subspace} \cite{FS1}
associated to $L(\Lambda)$, $\Lambda$ a dominant weight, is defined as
$$W(\Lambda)=U(\widehat{\goth{n}_+}) \cdot v_{\Lambda} \subset L(\Lambda),$$
where $v_{\Lambda}$ is a highest weight vector (which is unique up to
a nonzero multiple).

To study the principal subspaces it will be convenient, as in
\cite{CLM}, to use the polynomial algebra
$$\mathcal{A}=\mathbb{C}[y_{-1}, y_{-2},\dots],$$
where $y_{-1}, y_{-2}, \ldots $ are (commuting, independent) formal
variables.  For an $\widehat{\goth{sl}(2)}$--module $M$, consider the
algebra map
\bea
\mathcal{A} & \longrightarrow & {\rm End} \ M \nn
y_{-j} & \mapsto & x_{\alpha}(-j)
\eea
($j>0$); this map is well defined because the operators
$x_{\alpha}(-j)$ commute.
For a dominant weight $\Lambda$, define the linear surjection
\bea \label{fLambda}
f_{\Lambda}: \mathcal A & \rightarrow & W(\Lambda) \nn
p(y_{-1},y_{-2},\dots)
& \mapsto & p(x_{\alpha}(-1),x_{\alpha}(-2),\ldots) \cdot
v_{\Lambda},
\eea
where $p$ is a polynomial and where as above, $v_{\Lambda}$ is a
highest weight vector, and set
$$\mathcal A_{\Lambda}= {\rm Ker}\ f_{\Lambda},$$
an ideal in $\mathcal{A}$.  Then we have
\be
\mathcal A/\mathcal A_{\Lambda}
\stackrel{\sim}{\longrightarrow} W(\Lambda).
\ee

{}From now on we fix a positive integer $k$ (which corresponds to the
integer $l-1$ in the Introduction).

Our main goal is to derive explicit formulas for certain graded
dimensions of principal subspaces via certain systems of
$q$--difference equations. To achieve this goal we shall need an
explicit description of the ideals $\mathcal{A}_{\Lambda}$. The
following result was (essentially) proved in \cite{FS1}--\cite{FS2}:

\begin{theorem} \label{relationk}
For every $i$ with $0 \leq i \leq k$,
\begin{equation} \label{equality}
\mathcal{A}_{(k-i)\Lambda_0+i \Lambda_1}=
\mathcal{A}y_{-1}^{k-i+1}+\mathcal{A}_{k
\Lambda_0}.
\end{equation}
\end{theorem}
{\em Proof:} {}From \cite{FS1} it follows that the ideal
$$\mathcal{A}_{(k-i)\Lambda_0+i \Lambda_1}$$
is generated by the elements
$$r^{(k)}_n=\sum_{\stackrel{i_1,\dots,i_{k+1} > 0}{i_1+\cdots +
i_{k+1}=-n}} y_{-i_1} \cdots y_{-i_{k+1}}$$
for $n \leq -(k+1)$, and
$$y_{-1}^{k+1-i}.$$
This fact immediately implies the statement.
\epf

\renewcommand{\theequation}{\thesection.\arabic{equation}}
\renewcommand{\thetheorem}{\thesection.\arabic{theorem}}
\setcounter{equation}{0}
\setcounter{theorem}{0}

\section{Intertwining operators for vertex operator algebras
and fusion rules}

In this section we will be using the theory of vertex operator
algebras and intertwining operators, as developed in \cite{FLM},
\cite{FHL} and \cite{DL}.  The reader unfamiliar with the theory of
modules and intertwining operators for vertex operator algebras may
consult our previous paper \cite{CLM}, where we recalled, and
motivated, the definition of intertwining operator (see \cite{FHL} for
more details).

It is well known that $L(k \Lambda_0)$ (see the previous section) has
a natural vertex operator algebra structure. In addition, all the
level $k$ standard modules are modules for this vertex operator
algebra $L(k \Lambda_0)$, and conversely, these are all the
irreducible $L(k \Lambda_0)$--modules up to equivalence (see
\cite{FZ}; cf. \cite{Li1}, \cite{DL}, \cite{MP} and \cite{LL}).  Let
$L(\Lambda)$ be one of these irreducible $L(k \Lambda_0)$--modules.
It is known (cf. \cite{FZ}) that $L(\Lambda)$ is graded with respect
to a standard action of the Virasoro algebra element $L(0)$ (not to be
confused with the trivial $\widehat{\goth{sl}(2)}$--module!) and that
it decomposes as

\bea \label{grading}
L(\Lambda) &=&\coprod_{s \geq 0} L(\Lambda)_{s+h_{\Lambda}},
\nonumber
\eea
where
\bea \label{grading2}
L(\Lambda)_{\lambda}&=& \{ v \in L(\Lambda)\;|\; L(0) \cdot v=\lambda v
\} \nonumber
\eea
is the {\em weight} space of $L(\Lambda)$ of {\it weight} $\lambda \in
\Bbb C$, and where
$$h_{\Lambda}=\frac{\langle \Lambda,\Lambda+ \alpha
\rangle}{2(k+2)}$$
(cf. \cite{K}, \cite{DL}, \cite{LL}).  In addition, $L(\Lambda)$ has
second, compatible, grading, by {\em charge}, given by the eigenvalues
of the operator $\alpha(0)/2=h(0)/2$.  It is important to notice that
the principal subspace $W(\Lambda)$ is doubly graded as well. We will
denote by $L(\Lambda)_{r+ \langle \alpha/2,\Lambda \rangle,\; s +
h_{\lambda}}$ the subspace of $L(\Lambda)$ consisting of the vectors
of charge $r+ \langle \alpha/2,\Lambda \rangle$ and weight
$s+h_{\Lambda}$, so that
$$L(\Lambda)=\coprod_{r \in \mathbb{Z},\;s \in \mathbb{N}}
L(\Lambda)_{r+ \langle \alpha/2,\Lambda \rangle,\; s + h_{\Lambda}}.$$
Clearly, $W(\Lambda)$ also admits a decomposition
$$W(\Lambda)=\coprod_{r,s \in \mathbb{N}} W(\Lambda)_{r+ \langle
\alpha/2,\Lambda \rangle,\; s + h_{\Lambda}}.$$

Now let $V$ be an arbitrary vertex operator algebra
and let $W_1$, $W_2$ and $W_3$ be
$V$--modules.  We denote by
$$I \ { W_3 \choose W_1 \ W_2 }$$
the vector space of all intertwining operators of type
${ W_3 \choose W_1 \ W_2 }$ (see \cite{FHL}).
In many cases these spaces are finite--dimensional.  Their dimensions
$$\mathcal{N}_{W_1 \ W_2}^{W_3}= {\rm dim} \ I \ { W_3 \choose W_1 \
W_2 }$$
are the so--called {\em fusion coefficients} or {\em fusion rules}.
A formula for the fusion rules for irreducible modules for
$L(k \Lambda_0)$ is due to Frenkel and Zhu \cite{FZ} (for
clarification and generalization of Frenkel--Zhu's formula see
\cite{Li2}).  We should mention that the same result was obtained
previously in \cite{TK}, but not in the setting of vertex operator
algebras.  The following result is {}from \cite{FZ}:
\begin{proposition} \label{intpropo}
Let
$$\mathcal{N}_{i_1,i_2}^{i_3}={\rm dim} \ I { L((k-i_3)\Lambda_0+i_3
\Lambda_1) \choose L((k-i_1)\Lambda_0+i_1 \Lambda_1) \
L((k-i_2)\Lambda_0+i_2 \Lambda_1)} ,$$
where $0 \leq i_j \leq k$, $i_j \in \mathbb{N}$, $j=1,2,3.$
Then
$$\mathcal{N}_{i_1,i_2}^{i_3}=1$$
if and only if
$$i_1+i_2+i_3 \in 2\mathbb{Z},\;\; i_1+i_2+i_3 \leq 2k,\;\;
|i_1-i_2| \leq i_3 \leq i_1+i_2.$$
Otherwise, $\mathcal{N}_{i_1,i_2}^{i_3}=0.$
\end{proposition}

For example, if $k=1$, the nontrivial intertwining operators
are of the following types:
$${ L(\Lambda_0) \choose L(\Lambda_0) \ L(\Lambda_0) },\;
{ L(\Lambda_1) \choose L(\Lambda_0) \ L(\Lambda_1) },\;
{ L(\Lambda_0) \choose L(\Lambda_1) \ L(\Lambda_1) }
\ \mbox{and} \
{ L(\Lambda_1) \choose L(\Lambda_1) \ L(\Lambda_0) }.$$

In order to proceed we will need a more detailed description of the
intertwining operators indicated in Proposition \ref{intpropo}.
Let $W_1$, $W_2$ and $W_3$ be
irreducible $L(k \Lambda_0)$--modules.  Then it is not hard to see
(cf. \cite{FHL}, \cite{FZ}) that every
intertwining operator $\mathcal{Y}( \cdot, x)$ of type
$${W_3 \choose W_1 \ W_2}$$
satisfies the condition
$$\mathcal{Y}(w_1,x) \in x^{h_{W_3}-h_{W_1}-h_{W_2}}
{\rm End}(W_2,W_3)[[x,x^{-1}]]$$
for every $w_1 \in W_1$.
We will write (cf. \cite{FZ})
$$\mathcal{Y}(w_1,x)= x^{h_{W_3}-h_{W_1}-h_{W_2}} \sum_{n \in
\mathbb{Z}}
(w_1)_{[n]} x^{-n-1},$$
$$ (w_1)_{[n]} \in {\rm End}(W_2,W_3).$$
Also, for every $w_1 \in
W_1$, let
$$o_{\mathcal{Y}}(w_1)={\rm Coeff}_{x^0}x^{-h_{W_3}+h_{W_1}+h_{W_2}}
\mathcal{Y}(w_1,x)=(w_1)_{[-1]}.$$
(Note: If $h_{W_1}+h_{W_2}=h_{W_3}$, such as when $W_2=L(\Lambda_0)$
and $W_1=W_3$, then $o_{\mathcal{Y}}(w_1)$ is simply the constant term
of the intertwining operator $\mathcal{Y}(w_1,x)$.)
Let $v_{W_i}$ be a
highest weight vector of $W_i$, $i=1,2,3$.
Then for every nonzero $\mathcal{Y} \in I { W_3
\choose W_1 \ W_2}$,
\be \label{zeroth} o_{\mathcal{Y}}(v_{W_1})v_{W_2}=
(v_{W_1})_{[-1]}v_{W_2}, \ee
a highest weight vector of $W_3$ or zero; it will be nonzero for our
cases below. Also, {}from the commutator formula for intertwining
operators \cite{FHL}, it follows that
\be \label{comminter}
[x_{\alpha}(n),\mathcal{Y}(v_{W_1},x)]=0,
\ee
for
every $n \in \mathbb{Z}$. Therefore,
\be \label{commutingwithUn+^}
\label{comny} [U(\widehat{\goth{n}_+}),\mathcal{Y}(v_{W_1},x)]=0.
\ee

So far we have treated irreducible $L(k \Lambda_0)$--modules
modules abstractly. We have already outlined (cf. Section 2)
that the level $1$ standard modules admit an explicit
construction as subspaces of $V_P$
(where $P=\frac{1}{2}\mathbb{Z}\alpha$ is the weight lattice of
$\goth{sl}(2)$).
We can realize the level $k$ standard modules as submodules of
the tensor products of $k$ level $1$ standard modules
(cf. Section 2).  Let
$$\stackrel{\displaystyle{L=\underbrace{P \oplus \cdots \oplus P}} \; .}
{\ \ \ \ k \; \mbox{\small{times}} }$$
Consider the space
$V_P=V_Q \oplus V_{Q+{\alpha}/2}$ (cf.  \cite{FLM}, \cite{CLM}).
Let
$$\stackrel{\displaystyle{V_L=\underbrace{V_P \otimes \cdots \otimes
V_P}} \; ;}{\ \ \ \ \ \ k \; \mbox{\small{times}}  }$$
$V_L$ is naturally an $\widehat{\goth{sl}(2)}$--module.
Now, inside $V_P \cong L(\Lambda_0) \oplus L(\Lambda_1)$, we make the
following identifications:
$$1=v_{\Lambda_0} \ \ \mbox{and} \ \ e^{\alpha/2}=v_{\Lambda_1}.$$
There are of course many ways to embed
$L(k_0 \Lambda_0+k_1 \Lambda_1)$ ($k_0 + k_1 = k$) inside $V_L$.
Let $(i_1,\dots,i_k)$ be a $k$--tuple such that $i_j \in \{0,1\}$,
$j=1,\dots,k$.  Consider a vector of the form
$$v_{i_1,\dots,i_k}=
v_{\Lambda_{i_1}} \otimes \cdots \otimes v_{\Lambda_{i_k}} \in V_L,$$
where exactly $k_0$ indices (i.e., $i_j$'s)
are equal to $0$ (and exactly $k_1$ are equal
to $1$). This is certainly a highest weight vector for
$\widehat{\goth{sl}(2)}$, and
$$L(k_0 \Lambda_0 + k_1 \Lambda_1) \cong
U(\widehat{\goth{sl}(2)}) \cdot v_{i_1,\dots,i_k} \subset V_L$$
(cf. \cite{K}).
Let us denote by $\iota_{i_1,\dots,i_k}$
the embedding
$$\iota_{i_1,\dots,i_k} : L(k_0 \Lambda_0 + k_1 \Lambda_1)
\longrightarrow V_L,$$
uniquely determined by the identification
$$v_{k_0 \Lambda_0+k_1\Lambda_1}=v_{i_1,\dots,i_k}.$$

In parallel with \cite{CLM}, for $\lambda \in P$, we will consider
$$e^{\lambda}_{(k)} : V_L \longrightarrow V_L,$$
$$\stackrel{\displaystyle{e^{\lambda}_{(k)}=\underbrace{e^{\lambda}
\otimes \cdots \otimes e^{\lambda} } \; , } }{\ \ \ \ \ \
k \; \mbox{\small{times}}  }$$
a linear automorphism of $V_L$, with inverse $e^{-\lambda}_{(k)}$.
Then it follows (cf. \cite{Ge}, \cite{CLM}) that
$$
e^{\lambda}_{(k)} x_{\alpha}(-m_1) \cdots x_{\alpha}(-m_s) \cdot v=
x_{\alpha}(-m_1-\langle \lambda,\alpha \rangle ) \cdots
x_{\alpha}(-m_s-\langle \lambda,\alpha \rangle )
e^{\lambda}_{(k)} \cdot v,$$
for every $v \in V_L$ and $m_1,\dots,m_s \in \mathbb{Z}$.
Also,
$$ e^{\lambda}_{(k)} x_{-\alpha}(-m_1) \cdots x_{-\alpha}(-m_s) \cdot
v=x_{-\alpha}(-m_1+\langle \lambda,\alpha \rangle ) \cdots
x_{-\alpha}(-m_s+\langle \lambda,\alpha \rangle )
e^{\lambda}_{(k)} \cdot v$$
and if each $m_i \neq 0$,
$$ e^{\lambda}_{(k)} h(-m_1) \cdots h(-m_s) \cdot
v=h(-m_1) \cdots h(-m_s) e^{\lambda}_{(k)} \cdot v.$$

The following lemma describes the action of $e^{\alpha/2}_{(k)}$
on certain irreducible $\widehat{\goth{sl}(2)}$--submodules
of $V_L$ (i.e., irreducible $L(k \Lambda_0)$--modules)
and the corresponding principal subspaces:
\begin{lemma} \label{add}
\begin{itemize}
\item[(a)]
Consider the subspace
$L(k_0 \Lambda_0+k_1\Lambda_1) \subset V_L$ embedded via
$\iota_{i_1,\dots,i_k}$ as above.  The image of the restriction
map $e^{\alpha/2}_{(k)}|_{L(k_0 \Lambda_0+k_1\Lambda_1)}$
lies in
$L(k_1 \Lambda_0 + k_0 \Lambda_1)$, embedded via
$\iota_{1-i_1,\dots,1-i_k}$.
In particular, $e^{\alpha/2}_{(k)}$ defines maps
interchanging $L(k_0 \Lambda_0+ k_1 \Lambda_1)$ and
$L(k_1 \Lambda_0+k_0 \Lambda_1)$.
\item[(b)]
Consider the principal subspace
$W(k_0 \Lambda_0 + k_1 \Lambda_1) \subset L(k_0
\Lambda_0+k_1 \Lambda_1)$.
The image of the restriction
map $e^{\alpha/2}_{(k)}|_{W(k_0 \Lambda_0 +
k_1 \Lambda_1)}$ lies in $W(k_1 \Lambda_0 + k_0 \Lambda_1) \subset
L(k_1 \Lambda_0+k_0 \Lambda_1)$.
\item[(c)]
If $k_0=k$ (or $k_1=0$), the map in (b) is surjective and in
particular is a linear isomorphism {}from $W(k \Lambda_0)$ to
$W(k \Lambda_1)$.
\end{itemize}
\end{lemma}
{\em Proof:}
First we prove (a).  {}From the formulas preceding the lemma it
follows that
$$ e^{\alpha/2}_{(k)} \pi(U(\widehat{\goth{sl}(2)})) e^{-\alpha/2}_{(k)}=
\pi(U(\widehat{\goth{sl}(2)})),$$
where $\pi$ is our representation of $U(\widehat{\goth{sl}(2)})$ on $V_L$.
Therefore, we need only
show that $$e^{\alpha/2}_{(k)} \cdot v_{i_1,\dots,i_k} \in
U(\widehat{\goth{sl}(2)}) \cdot v_{1-i_1,\dots,1-i_k}.$$
For convenience let us assume that $i_1=i_2= \cdots =i_{k_1}=1$.
By using our identifications we have
$$e^{\alpha/2}v_{\Lambda_1}=e^\alpha \; (=e^{\alpha}v_{\Lambda_0})
=x_{\alpha}(-1)v_{\Lambda_0} \in V_P$$
(cf. \cite{CLM}).  Then
\bea \label{AA1}
&& \stackrel{e^{\alpha/2}_{(k)} \underbrace{v_{\Lambda_1} \otimes
\cdots \otimes v_{\Lambda_1} }
{\displaystyle {\otimes}}
\underbrace{v_{\Lambda_0} \otimes \cdots \otimes v_{\Lambda_0}
} \; =}{\ \ \ \ k_1 \; \mbox{\small{times}} \ \ \ \ \ \ \ \ \ \ \ \
k_0 \; \mbox{\small{times}} } 
\nonumber \\
&& \stackrel{ = \; \underbrace{ x_{\alpha}(-1)v_{\Lambda_0} \otimes
\cdots \otimes x_{\alpha}(-1) v_{\Lambda_0} }
{\displaystyle {\otimes}}
\underbrace{v_{\Lambda_1} \otimes \cdots \otimes v_{\Lambda_1}
}.}
{\ \ \ \ \ \ \ \ \ \ \ k_1 \; \mbox{\small{times}} \ \ \ \ \ \ \ \ \ \ 
\ \ \ \ \ \ \ \ \ \ \ \ k_0 \; \mbox{\small{times}}  }
\nonumber
\eea
Now, because
$$x_{\alpha}(-1)v_{\Lambda_1}=0 \ \ \mbox{and} \ \
x_{\alpha}^2(-1)v_{\Lambda_0}=0$$
(again cf. \cite{CLM}), it follows that
\begin{eqnarray}
&& \stackrel{\underbrace{x_{\alpha}(-1)v_{\Lambda_0} \otimes
\cdots \otimes x_{\alpha}(-1) v_{\Lambda_0} }
{\displaystyle {\otimes}}
\underbrace{v_{\Lambda_1} \otimes \cdots \otimes v_{\Lambda_1}} = }
{\ \ \ \ \ \ k_1 \; \mbox{\small{times}} \ \ \ \ \ \ \ \ \
\ \ \ \ \ \ \ \ \ \ \ \ \  k_0 \; \mbox{\small{times}} }
\nonumber \\
&&\stackrel{ \displaystyle{ =
{\frac{1}{{k_1}!}}\Delta(x_{\alpha}(-1)^{k_1})}
(\underbrace{v_{\Lambda_0}
\otimes
\cdots \otimes v_{\Lambda_0} }
\otimes \underbrace{v_{\Lambda_1} \otimes \cdots \otimes
v_{\Lambda_1}}),}
{\ \ \ \ \ \ \ \ \ \ \ \ \ \ \ \ \ \ \ \ \ \ k_1 \;
\mbox{\small{times}} \ \ \ \ \ \ \ \ \ \ \ \ \
k_0 \; \mbox{\small{times}}} \nonumber
\end{eqnarray}
so that
$$ e^{\alpha/2}_{(k)} \cdot v_{i_1,\dots,i_k} \in
U(\widehat{\goth{sl}(2)}) \cdot v_{1-i_1,\dots,1-i_k},$$
as desired.

The same reasoning with $U(\widehat{\goth{sl}(2)})$ replaced by
$U(\widehat{\goth{n}_+})$ proves (b).

To prove (c) it is enough to notice that
$$\stackrel{{e^{\alpha/2}_{(k)}} \; \cdot \; \underbrace{v_{\Lambda_0}
\otimes \cdots \otimes v_{\Lambda_0}} \; }{\ \ \ \ \ \ k \;
\mbox{\small{times}} }
\stackrel{= \ \underbrace{v_{\Lambda_1} \otimes \cdots \otimes
v_{\Lambda_1}} \; . \ \ \  \epf }{ k \; \mbox{\small{times}} }
$$

\section{The main theorem and consequences}
We continue to fix a positive integer $k$. Let $\Lambda$ be a dominant
weight such that $\langle \Lambda,c \rangle=k$. We shall determine the
following graded dimensions (generating functions of the dimensions of
homogeneous subspaces) of the spaces $W(\Lambda)$:
$$\chi_{W(\Lambda)}(x,q)={\rm dim}_*(W(\Lambda),x,q)={\rm
tr}|_{W(\Lambda)} x^{\alpha/2} q^{L(0)},$$
where we are using the formal variables $x$ and $q$.  To avoid the
multiplicative factor $x^{\langle \alpha/2, \Lambda
\rangle}q^{h_{\Lambda}}$ it is convenient to use slightly modified
graded dimensions: We define
$${\chi}'_{W(\Lambda)}(x,q)=x^{-\langle \alpha/2, \Lambda
\rangle}q^{-h_{\Lambda}} {\chi}_{W(\Lambda)}(x,q),$$
so that
$${\chi}'_{W(\Lambda)}(x,q) \in {\Bbb C}[[x,q]],$$
and in fact,
$${\chi}'_{W(\Lambda)}(x,q) \in 1 + xq{\Bbb C}[[x,q]].$$
It is a nontrivial task to find explicit formulas for the
${\chi}'_{W(\Lambda)}(x,q)$. As in \cite{CLM}, we will first derive
certain representation--theoretic results and convert these results
into statements about graded traces and $q$--difference equations.

The simplest $q$--difference equation follows easily {}from Lemma
\ref{add}, part (c), just as in formula (3.26) of \cite{CLM}:

\begin{proposition} \label{i0lemma}
We have
$${\chi}_{W(k \Lambda_1)}(x,q)= x^{k/2} q^{k/4} {\chi}_{W(k
\Lambda_0)}(xq,q)= x^{k/2} q^{h_{k \Lambda_1}} {\chi}_{W(k
\Lambda_0)}(xq,q),$$
or simply
$${\chi}'_{W(k \Lambda_1)}(x,q)={\chi}'_{W(k\Lambda_0)}(xq,q).
\ \ \  \epf $$
\end{proposition}

Our main result is:

\begin{theorem} \label{main}
Let $1 \leq i \leq k$, let
$$\mathcal{Y}( \cdot, x) \  \in I \ {L((i-1)\Lambda_0+(k-i+1)\Lambda_1)
\choose
L((k-1)\Lambda_0+\Lambda_1) \ \ L(i
\Lambda_0+(k-i)\Lambda_1) }$$
be a nonzero intertwining operator and let
$o_{\mathcal{Y}}(v_{(k-1)\Lambda_0+\Lambda_1})$
be as in (\ref{zeroth}).
Then the sequence
\begin{eqnarray} \label{complexk}
&& 0 \longrightarrow W((k-i)\Lambda_0+i \Lambda_1)
\stackrel{e^{\alpha/2}_{(k)}}{\longrightarrow} W(i
\Lambda_0+(k-i)\Lambda_1)
\stackrel{o_{\mathcal{Y}}(v_{(k-1)\Lambda_0+\Lambda_1})}{\longrightarrow}
\nonumber \\
&& \ \ \ \ \ \longrightarrow W((i-1)\Lambda_0+(k-i+1)\Lambda_1)
\longrightarrow 0
\end{eqnarray}
is exact. Also,
\begin{equation} \label{exact0}
0 \longrightarrow W(k\Lambda_0)
\stackrel{e^{\alpha/2}_{(k)}}{\longrightarrow}
W(k\Lambda_1) \longrightarrow 0
\end{equation}
is exact.
\end{theorem}

\begin{remark} \label{i0}
{\em When $k=i=1$, this theorem is equivalent to the main theorem in
\cite{CLM}, which yielded the Rogers--Ramanujan recursion.}
\end{remark}

\noindent 
{\em Proof of Theorem \ref{main}:} The exactness of (\ref{exact0}) has
been proved in Lemma \ref{add}. Now we prove the exactness of
(\ref{complexk}). First note that a nonzero intertwining operator
$\mathcal{Y}( \cdot, x)$ exists by Proposition \ref{intpropo}.  We
recall (see (\ref{zeroth})) that
$$o_{\mathcal{Y}}(v_{(k-1)\Lambda_0+\Lambda_1})v_{i
\Lambda_0+(k-i)\Lambda_1}$$
is a nonzero multiple of $v_{(i-1)\Lambda_0+(k-i+1)\Lambda_1}.$ In
addition, $o_{\mathcal{Y}}(v_{(k-1)\Lambda_0+\Lambda_1})$ commutes
with the action of ${U}(\widehat{\goth{n}_+})$; thus
$o_{\mathcal{Y}}(v_{(k-1)\Lambda_0+\Lambda_1})$ is surjective.  We
already know that $e^{\alpha/2}_{(k)}$ is injective and that it maps
$W((k-i)\Lambda_0+i \Lambda_1)$ into $W(i\Lambda_0+(k-i)\Lambda_1)$
(Lemma \ref{add}).

Let us prove the chain property.  First, we have
\begin{equation} \label{ex}
e^{\alpha/2}_{(k)} x_{\alpha}(-m_1) \cdots x_{\alpha}(-m_n)
=x_{\alpha}(-m_1-1) \cdots x_{\alpha}(-m_n-1)e^{\alpha/2}_{(k)}
\end{equation}
for $m_j \in \mathbb{Z}$.
By combining (\ref{ex}) with the proof of Lemma \ref{add}, we see that
the image of $e^{\alpha/2}_{(k)}$ is the
$U(\widehat{\goth{n}_+})$--submodule of $W(i\Lambda_0+(k-i)\Lambda_1)$
generated by
$$
\stackrel{ \displaystyle{ \Delta(x_{\alpha}(-1)^i)}
(\underbrace{v_{\Lambda_0} \otimes \cdots \otimes v_{\Lambda_0} }
\otimes \underbrace{v_{\Lambda_1} \otimes \cdots \otimes
v_{\Lambda_1} }).} {\ \ \ \ \ \ \ \ \ \ \ \ \ \  \ i \
\mbox{\small{times}} \ \ \ \ \ \ \ \ \ \ \ \ k-i \
\mbox{\small{times}} }
$$
But by (\ref{comminter}) (or (\ref{commutingwithUn+^})) and
(\ref{zeroth}),
\begin{eqnarray}
&& \stackrel{\displaystyle{o_{\mathcal{Y}}(v_{(k-1)\Lambda_0+\Lambda_1})
\Delta(x_{\alpha}(-1)^i)}
(\underbrace{v_{\Lambda_0} \otimes \cdots \otimes v_{\Lambda_0} }
\otimes \underbrace{v_{\Lambda_1} \otimes \cdots \otimes
v_{\Lambda_1} })} {\ \ \ \ \ \ \ \ \ \ \ \ \ \ \ \ \ \ \ 
\ \ \ \ \ \ \ \ \ \ \ \ \ \ \ \ \ \ \ i \
\mbox{\small{times}} \ \ \ \ \ \ \ \ \ \ \ \ k-i \
\mbox{\small{times}} } 
\nonumber \\
&& \stackrel{\displaystyle{= \Delta(x_{\alpha}(-1)^i)}
o_{\mathcal{Y}}(v_{(k-1)\Lambda_0+\Lambda_1})
(\underbrace{v_{\Lambda_0} \otimes \cdots \otimes v_{\Lambda_0} }
\otimes \underbrace{v_{\Lambda_1} \otimes \cdots \otimes
v_{\Lambda_1} })} {\ \ \ \ \ \ \ \ \ \ \ \ \ \ \ \ \ \ \ \ 
\ \ \ \ \ \ \ \ \ \ \ \ \ \ \ \ \ \ \ \ i \
\mbox{\small{times}} \ \ \ \ \ \ \ \ \ \ \ \ k-i \
\mbox{\small{times}} } 
\nonumber \\
&& \stackrel{\displaystyle{=a \Delta(x_{\alpha}(-1)^i)}
(\underbrace{v_{\Lambda_0} \otimes \cdots \otimes v_{\Lambda_0} }
\otimes \underbrace{v_{\Lambda_1} \otimes \cdots \otimes
v_{\Lambda_1} }),} {\ \ \ \ \ \ \ \ \ \ \ \ \ \ \ \ \ \ \ \ i-1 \
\mbox{\small{times}} \ \ \ \ \ \ \ k-i+1 \
\mbox{\small{times}} }
\nonumber
\end{eqnarray}
where $a$ is a nonzero constant, and this equals $0$ in view of (the
easy part of) Theorem \ref{relationk}.  The chain property now follows
by another use of (\ref{comminter}).

Finally, to prove the exactness, we continue to follow \cite{CLM}.
First we characterize the kernel of the map
$o_{\mathcal{Y}}(v_{(k-1)\Lambda_0+\Lambda_1})$.  Suppose that $p$ is
a polynomial in $\mathcal{A}$, so that $f_{i
\Lambda_0+(k-i)\Lambda_1}(p)$ is a general element of $W(i
\Lambda_0+(k-i)\Lambda_1)$ (recall (\ref{fLambda})).  Then
$$o_{\mathcal{Y}}(v_{(k-1)\Lambda_0+\Lambda_1})
(f_{i \Lambda_0+(k-i)\Lambda_1}(p))=0$$
if and only if
$$p \in {\rm Ker} \ f_{(i-1)\Lambda_0+(k-i+1)\Lambda_1}=
\mathcal{A}_{(i-1)\Lambda_0+(k-i+1)\Lambda_1}.$$
Now, the description of the ideals in (\ref{equality}) implies that
$$\mathcal{A}_{ (i-1)\Lambda_0+(k-i+1)\Lambda_1 }
=\mathbb{C}[y_{-2},y_{-3},\dots]y_{-1}^i + \mathcal{A}_{ i
\Lambda_0+(k-i)\Lambda_1}.$$
Thus
\be \label{left}
o_{\mathcal{Y}}(v_{(k-1)\Lambda_0+\Lambda_1})
(f_{i \Lambda_0+(k-i)\Lambda_1}(p))=0
\Longleftrightarrow
p \in \mathbb{C}[y_{-2},y_{-3},\dots]y_{-1}^i +
\mathcal{A}_{i \Lambda_0+(k-i)\Lambda_1}.
\ee

Next we characterize the image of
$$e^{\alpha/2}_{(k)} : W((k-i)\Lambda_0+i \Lambda_1)
\longrightarrow
W(i \Lambda_0+(k-i)\Lambda_1).$$
Notice that $f_{i \Lambda_0+(k-i)\Lambda_1}(p)$ is
of the form  $e^{\alpha/2}_{(k)}(w)$ for some $w \in
W((k-i)\Lambda_0+i \Lambda_1)$
if and only if for some
$q=q(y_{-1},y_{-2},\dots) \in \mathcal{A}$,
$$f_{i \Lambda_0+(k-i)\Lambda_1}(p)=
e^{\alpha/2}_{(k)}(f_{(k-i)\Lambda_0+i \Lambda_1}(q))$$
$$=q(x_{\alpha}(-2),x_{\alpha}(-3),\dots)e^{\alpha/2}_{(k)}
v_{(k-i)\Lambda_0+i \Lambda_1}$$
$$={\frac{1}{i!}}
q(x_{\alpha}(-2),x_{\alpha}(-3),\dots)x_{\alpha}(-1)^i
\cdot v_{i \Lambda_0+(k-i)\Lambda_1}$$
$$={\frac{1}{i!}}
f_{i\Lambda_0+(k-i) \Lambda_1}
(q(y_{-2},y_{-3},\dots)y_{-1}^i).$$
But this holds if and only if
$$p-q(y_{-2},y_{-3},\dots)y_{-1}^i \in
\mathcal{A}_{i\Lambda_0+(k-i) \Lambda_1}.$$
In other words,
\be \label{right}
f_{(k-i)\Lambda_0+i \Lambda_1}(p) \in {\rm Im}(
e^{\alpha/2}_{(k)}) \Longleftrightarrow
p \in \mathbb{C}[y_{-2},y_{-3},\dots]y_{-1}^i+
\mathcal{A}_{i \Lambda_0+(k-i)\Lambda_1}.
\ee
The exactness now follows {}from
(\ref{left}) and (\ref{right}).
\epf

Let us derive a few consequences of the theorem.  As above, we denote
by
$$W(\Lambda)_{r+\langle \alpha/2,\Lambda \rangle,\;
s+h_{\Lambda} } \subset W(\Lambda)$$
the subspace of $W(\Lambda)$ consisting of the vectors of charge $r
+ \langle \alpha/2, \Lambda \rangle$ and of weight
$s+h_{\Lambda}$.  Write
$$W'(\Lambda)_{r,s}=W(\Lambda)_{r+\langle \alpha/2,\Lambda
\rangle,\; s+h_{\Lambda}}.$$
{}From our construction it is clear that for $i=0,\dots,k,$
$e^{\alpha/2}_{(k)}$ maps $W'((k-i)\Lambda_0+i \Lambda_1)_{r,s}$ to
$W'(i \Lambda_0+ (k-i) \Lambda_1)_{r+i,s+r+i}$. Also, for
$i=1,\dots,k,$ $o_{\mathcal{Y}}(v_{(k-1)\Lambda_0+\Lambda_1})$ maps
$W'(i \Lambda_0+ (k-i) \Lambda_1)_{r+i,s+r+i}$ to $W'((i-1) \Lambda_0+
(k-i+1) \Lambda_1)_{r+i,s+r+i}$. Now, because of the exactness,
for $i=1,\dots,k,$
$$\frac{W'(i \Lambda_0+ (k-i)
\Lambda_1)_{r+i,s+r+i}}{e^{\alpha/2}_{(k)}
W'((k-i)\Lambda_0+i \Lambda_1)_{r,s}}
\cong W'((i-1) \Lambda_0+ (k-i+1) \Lambda_1)_{r+i,s+r+i}$$
(and similarly for $i=0$), and so
\bea
&& {\rm dim} \ \{ W'(i \Lambda_0+ (k-i) \Lambda_1)_{r+i,s+r+i} \} -
{\rm dim} \ \{W'((k-i)\Lambda_0+i \Lambda_1)_{r,s} \}
\nonumber \\
&& ={\rm dim} \ \{ W'((i-1) \Lambda_0+ (k-i+1) \Lambda_1)_{r+i,s+r+i}
\}, \nonumber 
\eea
for every $r$ and $s$. The last formula written in generating function
form yields 
$$x^{-i}{\chi}'_{W(i\Lambda_0+(k-i)\Lambda_1)}(x/q,q)-{\chi}'_{W((k-i)
\Lambda_0+ i\Lambda_1)}(x,q)$$
$$=x^{-i}{\chi}'_{W((i-1)\Lambda_0+(k-i+1)\Lambda_1)}(x/q,q).$$
If we multiply this equation by $x^i$ and substitute $xq$ for $x$, we
obtain
\begin{eqnarray} \label{recursions1}
&& {\chi}'_{W(i\Lambda_0+(k-i)\Lambda_1)}(x,q)-(xq)^i
{\chi}'_{W((k-i)\Lambda_0+i\Lambda_1)}(xq,q) \nn
&& = {\chi}'_{W((i-1)\Lambda_0+ (k-i+1)\Lambda_1)}(x,q),
\end{eqnarray}
for every $1 \leq i \leq k$,
and
\begin{eqnarray} \label{recursions2}
{\chi}'_{W(k\Lambda_1)}(x,q)={\chi}'_{W(k\Lambda_0)}(xq,q).
\end{eqnarray}
These are exactly the Rogers--Selberg recursion formulas mentioned in
the Introduction; cf. \cite{A1} and Lemma 7.2 of \cite{A}, with our
${\chi}'_{W(i\Lambda_0+(k-i)\Lambda_1)}(x,q)$ playing the role of
$J_{k+1,i+1}(0;x;q)$ in that Lemma.  Of course, in these earlier
works, these recursions applied to series in $x$ and $q$, not to the
graded dimensions of graded vector spaces.

\begin{remark}
{\em Historically, recursion formulas equivalent to
(\ref{recursions1}) appeared for the first time in \cite{RR}.  They
were independently rediscovered by Selberg \cite{S}.  In \cite{A1}
Andrews used these recursions to give a new proof of the Gordon
identities, and in \cite{A4}, to prove his ``analytic form'' (multisum
form) of the Gordon identities; cf. \cite{A}.}
\end{remark}

It is easy to see, as in \cite{A}, that the recursion relations
(\ref{recursions1}) and (\ref{recursions2}), together with the initial
conditions $\lim_{x \to 0} \chi'_{W(\Lambda)}(x,q)=1$ and $\lim_{q \to
0} \chi'_{W(\Lambda)}(x,q)=1$ (here the limits stand for the formal
substitutions), uniquely determine the graded dimensions
${\chi}'_{W(i\Lambda_0+(k-i)\Lambda_1)}(x,q)$, $i=0,\dots,k$, as
elements of $\mathbb{C}[[x,q]]$.

\begin{example}
{\em Let $k=2$. We have three $q$--difference equations:
\bea
&&
{\chi}'_{W(\Lambda_0+\Lambda_1)}(x,q)- xq
{\chi}'_{W(\Lambda_0+\Lambda_1)}(xq,q)={\chi}'_{ W(2 \Lambda_1)}(x,q), \nn
&&
{\chi}'_{W(2\Lambda_0)}(x,q)-(xq)^2{\chi}'_{W(2\Lambda_1)}(xq,q)={\chi}'_{
W(\Lambda_0+\Lambda_1)}(x,q) \nonumber
\eea
and}
\begin{eqnarray}
{\chi}'_{W(2\Lambda_0)}(xq,q)&=&{\chi}'_{W(2\Lambda_1)}(x,q). 
\nonumber
\end{eqnarray}
\end{example}

We finally have:

\begin{corollary}
For every $i = 0,\dots,k$,
\bea \label{analytic}
{\chi}'_{W(i \Lambda_0+(k-i)\Lambda_1)}(x,q)=
\sum_{m \geq 0} \sum_{\stackrel{N_1+\cdots +N_{k}=m}{N_1 \geq \cdots \geq N_{k}
\geq 0}}
\frac{x^{m} q^{N_1^2 + \cdots + N_{k}^2+N_{i+1}+\cdots +
N_{k}}}{(q)_{N_1-N_2} \cdots (q)_{N_{k-1}-N_{k}}
(q)_{N_{k}}},
\eea
or
\bea
&& {\chi}_{W(i \Lambda_0+(k-i)\Lambda_1)}(x,q)= \nonumber \\
&& \sum_{m \geq 0} \sum_{\stackrel{N_1+\cdots +N_{k}=m}{N_1 \geq \cdots \geq N_{k}
\geq 0}}
\frac{x^{m+(k-i)/2} q^{h_{i \Lambda_0+(k-i)\Lambda_1}+N_1^2 + \cdots +
N_{k}^2+N_{i+1}+\cdots + N_{k}}}{(q)_{N_1-N_2} \cdots (q)_{N_{k-1}-N_{k}}
(q)_{N_{k}}}. \nonumber
\eea
\end{corollary}
{\em Proof:}
As in \cite{A4}, \cite{A}, the expressions on the right--hand side of
(\ref{analytic}) satisfy the system of recursions (\ref{recursions1}),
(\ref{recursions2}).  Because of the uniqueness the result follows.
\epf

The last corollary gives the Feigin--Stoyanovsky character formulas
obtained in \cite{FS1} (formula 2.3.3$'$) and in \cite{Ge}.

\noindent {\small \sc
Dipartimento Me. Mo. Mat.,
Universit{\`a} di Roma ``La Sapienza,''
Via A. Scarpa 16,
00161 Roma, Italia} \\
{\em E--mail address}: capparelli@dmmm.uniroma1.it \\

\vspace{2mm}
\noindent {\small \sc Department of Mathematics,  Rutgers University,
Piscataway,
NJ 08854} \\
{\em E--mail address}: lepowsky@math.rutgers.edu \\

\vspace{2mm}
\noindent {\small \sc Department of Mathematics,  University of Arizona,
Tucson, AZ 85721} \\
{\small Current address: \sc Department of Mathematics and Statistics,
University at Albany (SUNY), Albany, NY 12222} \\
{\em E--mail address}: amilas@math.albany.edu

\end{document}